\documentclass[11pt]{article}
\usepackage{graphicx}
\usepackage{epsfig}
\usepackage{latexsym}
\newsavebox{\toy}
\savebox{\toy}{\framebox[0.65em]{\rule{0cm}{1ex}}}
\newcommand{\QED}{\usebox{\toy}\end{demo}}
\newenvironment{property}%
{\begin{list}{}{\setlength{\rightmargin}{0pt}%
\setlength{\itemsep}{0pt}}}{\end{list}}
\newlength{\templength}
\newcommand{\bp}{\setlength{\templength}{\labelwidth}%
\setlength{\labelwidth}{2em}\begin{property}}
\newcommand{\ep}{\end{property}\setlength{\labelwidth}{\templength}}
\newtheorem{theorem}{Theorem}[subsection]
\newtheorem{lemma}[theorem]{Lemma}
\newtheorem{proposition}[theorem]{Proposition}
\newtheorem{corollary}[theorem]{Corollary}
\newtheorem{assumption}{Assumption}
\newtheorem{definition}{Definition}[subsection]
\newtheorem{remark}{Remark}[subsection]
\newtheorem{exercise}{Exercise}[subsection]

\newcommand{\Thm}[1]{Theorem \ref{Thm.#1}}
\newcommand{\Lem}[1]{Lemma \ref{Lem.#1}}

\newcommand{\Cor}[1]{Corollary \ref{Cor.#1}}

\newcommand{\Theorem}[1]{\begin{theorem}\label{Thm.#1}}
\newcommand{\Lemma}[1]{\begin{lemma}\label{Lem.#1}}
\newcommand{\Proposition}[1]{\begin{proposition}\label{Prop.#1}}
\newcommand{\Corollary}[1]{\begin{corollary}\label{Cor.#1}}
\newcommand{\Assumption}[1]{\begin{assumption}\label{Ass.#1}\rm}
\newcommand{\Definition}[1]{\begin{definition}\label{Def.#1}\rm}
\newcommand{\Remark}[1]{\begin{remark}\label{Rem.#1}\rm }
\newcommand{\Exercise}[1]{\begin{exercise}\label{Exe.#1}\rm }

\newcommand{\bd}{\begin{displaymath}}
\newcommand{\ed}{\end{displaymath}}
\newcommand{\bdn}{\begin{equation}}
\newcommand{\bdnl}{\begin{equation}\label}
\newcommand{\edn}{\end{equation}}
\newcommand{\barray}{\begin{array}}
\newcommand{\earray}{\end{array}}
\newcommand{\bds}{\begin{description}}
\newcommand{\eds}{\end{description}}
\newcommand{\bitemize}{\begin{itemize}}
\newcommand{\eitemize}{\end{itemize}}
\newcommand{\benumerate}{\begin{enumerate}}
\newcommand{\eenumerate}{\end{enumerate}}
\newcommand{\btabbing}{\begin{tabbing}}
\newcommand{\etabbing}{\end{tabbing}}
\newcommand{\bcenter}{\begin{center}}
\newcommand{\ecenter}{\end{center}}
\newcommand{\bflushright}{\begin{flushright}}
\newcommand{\bflushleft}{\begin{flushleft}}
\newcommand{\eflushright}{\end{flushright}}
\newcommand{\eflushleft}{\end{flushleft}}
\newcommand{\bdnn }{\begin{eqnarray*}}
\newcommand{\ednn }{\end{eqnarray*}}
\newcommand{\bdmn}{\begin{eqnarray}}
\newcommand{\edmn}{\end{eqnarray}}

\newcommand{\nn}{\nonumber}
\newcommand{\SSC}[1]{\section{#1}\setcounter{equation}{0}}

\newcounter{biblio}
\newenvironment{references}%
{\begin{list}{[\arabic{biblio}]}{\usecounter{biblio}%
\setlength{\leftmargin}{2.5em}\setlength{\rightmargin}{0pt}%
\setlength{\labelwidth}{2em}\setlength{\itemsep}{0pt}}}{\end{list}}
\newcommand{\References}%
{\vspace{2.8ex plus .3ex minus .3ex}%
\begin{center}{\bf References}\end{center}\begin{references}}

\usepackage{amssymb}
\usepackage{mathrsfs}
\newcommand{\bL}{{\mathbb{L}}}
\newcommand{\N}{{\mathbb{N}}}
\newcommand{\Z}{{\mathbb{Z}}}
\newcommand{\zd}{\Z^d}

\newcommand{\R}{{\mathbb{R}}}
\newcommand{\rd}{\R^d}

\newcommand{\ra }{\rightarrow }
\newcommand{\lra }{\longrightarrow }

\newcommand{\Ra}{\Rightarrow }

\newcommand{\LRa}{\Leftrightarrow }
\newcommand{\ov}{\overline}
\newcommand{\tl}{\widetilde}

\newcommand{\vvs}{\vspace{2ex}}
\newcommand{\vs}{\vspace{1ex}}

\newcommand{\lef}{\left}
\newcommand{\rig}{\right}
\newcommand{\ri}{\right}
\newcommand{\st}{\stackrel}

\newcommand{\8}{\infty}

\newcommand{\dps}{\displaystyle}
\newcommand{\sub}{\subset}

\newcommand{\bsh}{\backslash}

\newcommand{\half}{\mbox{$\frac{1}{2}$}}

\newcommand{\inflim}{\mathop{\underline{\lim}}}
\newcommand{\suplim}{\mathop{\overline{\lim}}}

\renewcommand{\a}{\alpha}

\newcommand{\Gm}{\Gamma}

\newcommand{\del}{\delta}

\newcommand{\e}{\varepsilon}

\newcommand{\h}{\eta}

\newcommand{\kp}{\kappa}

\newcommand{\lm}{\lambda}

\newcommand{\n}{\nu}

\newcommand{\rh}{\rho}

\newcommand{\s}{\sigma}

\renewcommand{\t}{\tau}

\newcommand{\W}{\Omega}

\newcommand{\cE }{{\cal E}}
\newcommand{\cF }{{\cal F}}

\newcommand{\cO }{{\cal O}}

\newcommand{\cR }{{\cal R}}

\makeatletter
\def\section{\@startsection{section}{1}{\z@}{-3.5ex plus -1ex minus 
 -.2ex}{2.3ex plus .2ex}{\bf}}
\def\subsection{\@startsection{subsection}{2}{\z@}{-3.25ex plus -1ex minus 
 -.2ex}{1.5ex plus .2ex}{\bf}}
\makeatother

\begin{document}
\bcenter

\large{\bf Central Limit Theorem for a Class 
of Linear Systems}\footnote{submitted on October 4, 2008. 
final version accepted on April 3, 2009.}

\vvs \normalsize

\vvs
Yukio Nagahata \\
Department of Mathematics, \\
Graduate School of Engineering Science \\
Osaka University,\\
Toyonaka 560-8531, Japan.\\
email: {\tt nagahata@sigmath.es.osaka-u.ac.jp}\\
URL: {\tt http://www.sigmath.osaka-u.ac.jp/}$\widetilde{}$ {\tt nagahata/}

\vvs
Nobuo Yoshida\footnote{
Supported in part by JSPS Grant-in-Aid for Scientific
Research, Kiban (C) 17540112} \\
Division of Mathematics \\
Graduate School of Science \\
Kyoto University,\\
Kyoto 606-8502, Japan.\\
email: {\tt nobuo@math.kyoto-u.ac.jp}\\
URL: {\tt http://www.math.kyoto-u.ac.jp/}$\widetilde{}$ {\tt nobuo/}

\ecenter

\begin{abstract}
We consider a class of interacting particle systems with values 
in $[0,\8)^{\zd}$, of which the binary contact path process 
is an example. 
For $d \ge 3$ and under a certain square integrability condition 
on the total number of the particles, we prove a central limit theorem 
for the density of the particles, together 
with upper bounds for the density of the most populated site and 
the replica overlap. 
\end{abstract}

\small
\noindent Abbreviated Title:  CLT for Linear Systems.
\tableofcontents

\normalsize
\SSC{Introduction}
We write $\N=\{0,1,2,...\}$, 
$\N^*=\{1,2,...\}$ and 
$\Z=\{ \pm x \; ; \; x \in \N \}$. For 
$x=(x_1,..,x_d) \in \rd$, $|x|$ stands for the $\ell^1$-norm: 
$|x|=\sum_{i=1}^d|x_i|$. For $\h=(\h_x)_{x \in \zd} \in \R^{\zd}$, 
$|\h |=\sum_{x \in \zd}|\h_x|$. 
Let 
$(\W, \cF, P)$ be 
a probability space.
We write $P[X]=\int X \; dP$ and 
$P[X:A]=\int_A X \; dP$ for a r.v.(random variable) $X$ and an 
event $A$. 
\subsection{The binary contact path process (BCPP)}
We start with a motivating simple example. 
Let $\h_t =(\h_{t,x})_{x \in \zd} \in \N^{\zd}$, $t \ge 0$ be 
binary contact path process (BCPP for short) with parameter $\lm>0$.
Roughly speaking, the BCPP is an extended version of the basic 
contact process, in which not only the presence/absence of 
the particles at each site, but also their number is considered. 
The BCPP was originally introduced by D. Griffeath 
\cite{Gri83}. Here, we explain the process following 
the formulation in the book of T. Liggett \cite[Chapter IX]{Lig85}.
Let $\t^{z,i}$, ($z \in \zd$, $i \in \N^*$) 
be i.i.d. mean-one exponential random variables 
and $T^{z,i}=\t^{z,1}+...+\t^{z,i}$. 
We suppose that the process $(\h_t)$ 
starts from a deterministic configuration 
$\h_0=(\h_{0,x})_{x \in \zd}\in \N^{\zd}$  
with $|\h_0|<\8$. 
At time 
$t=T^{z,i}$, $\h_{t-}$ is replaced by $\h_t$ randomly 
as follows: for each $e \in \zd$ with $|e|=1$, 
$$
\h_{t,x}=\lef\{ \barray{ll}
\h_{t-,x}+\h_{t-,z} & \mbox{if $x = z+e$,}\\
\h_{t-,x} & \mbox{if otherwise} \\
\earray \ri.\; \; 
\mbox{with probability ${\lm \over 2d\lm +1}$,} 
$$
(all the particles at site $z$ are 
duplicated and added to those on the site $z=x+e$), and 
$$
\h_{t,x}=\lef\{ \barray{ll}
0 & \mbox{if $x=z$}, \\
\h_{t-,x} & 
\mbox{if $x \neq z$}
\earray \ri.\; \; 
\mbox{with probability ${1 \over 2d\lm +1}$}
$$
(all the particles at site $z$ disappear).
The replacement occurs independently for different $(z,i)$ and independently 
from $\{ \t^{z,i}\}_{z,i}$. A motivation to study the BCPP comes from the 
fact that the projected process
$$
\lef( \h_{t,x} \wedge 1\ri)_{x \in \zd},\; \; \; t \ge 0
$$
is the basic contact process \cite{Gri83}.

Let 
$$
\kp_1={2d\lm -1 \over 2d\lm +1}\; \; \; \mbox{and}\; \; \;  
\ov{\h}_t=(\exp(-\kp_1t)\h_{t,x})_{x \in \zd}.
$$
Then,  $(|\ov{\h}_t|)_{t \ge 0}$
is a nonnegative martingale and therefore, the 
following limit exists almost surely:
$$
|\ov{\h}_\8|\st{\rm def}{=}
\lim_t|\ov{\h}_t|.
$$
Moreover, $P [|\ov{\h}_\8|]=1$ if 
\bdnl{|h_8|>0}
\mbox{$d \ge 3$ and $\lm > {1 \over 2d (1-2\pi_d)}$,} 
\edn
where $\pi_d$ is the return probability for the simple random 
walk on $\zd$ \cite[Theorem 1]{Gri83}. 
It is known that $\pi_d \le \pi_3=0.3405...$ for $d \ge 3$ 
\cite[page 103]{Spi76}.

We denote the
density of the particles by: 
\bdnl{rh}
\rh_{t,x}=\frac{\h_{t,x}}{|\h_t|}{\bf 1}\{ |\h_t| >0\}, \;
\; t >0, x \in \zd. 
\edn 
Interesting objects related to the
density would be 
\bdnl{rh^*} \rh^*_t=\max_{x \in \zd}\rh_{t,x}, \;
\; \mbox{and}\; \; \cR_t=\sum_{x \in \zd}\rh_{t,x}^{2}. 
\edn
$\rh^*_t$ is the density at the most populated site, while $\cR_t$
is the probability that a given pair of particles at time $t$ are
at the same site. We call $\cR_t$ the {\it replica overlap}, in
analogy with the spin glass theory. Clearly, $(\rh^*_t)^{2} \le \cR_t
\le \rh^*_t$. These quantities convey information on
localization/delocalization of the particles. Roughly speaking,
large values of $\rh^*_t$ or $\cR_t$ indicate that the most of
the particles are concentrated on small number of ``favorite
sites" ({\it localization}), whereas small values of them imply
that the particles are spread out over a large number of sites ({\it
delocalization}).

As  a special case of \Cor{CLT} below, we have 
the following result, which shows the diffusive behavior 
and the delocalization of the BCPP under the 
condition (\ref{|h_8|>0}):
\Theorem{BCPP}
Suppose (\ref{|h_8|>0}). Then, for any $f \in C_{\rm b}(\rd)$, 
$$
\lim_{t \ra \8} \sum_{x \in \zd}
f\lef(x/\sqrt{t}\ri)\rh_{t,x}
=\int_{\rd}fd\n\; \; \; 
\mbox{in $P(\; \cdot \; | |\ov{\h}_\8|>0)$-probability,}
$$
where 
$C_{\rm b} (\rd)$ 
stands for the set of bounded continuous functions on $\rd$, and 
$\n$ is the Gaussian measure with
$$
\int_{\rd}x_id\n (x)=0, 
\; \; \; \int_{\rd}x_ix_jd\n (x)={\lm \over 2d\lm +1}\del_{ij},\; \; \; 
i,j=1,..,d.
$$
Furthermore, 
$$
\cR_t=\cO (t^{-d/2})\; \; \; 
\mbox{as $t \nearrow \8$ in $P(\; \cdot \; | |\ov{\h}_\8|>0)$-probability}.
$$
\end{theorem}
\subsection{The results}
We generalize \Thm{BCPP} to a certain class of linear interacting 
particle systems with values in $[0,\8)^{\zd}$\cite[Chapter IX]{Lig85}.
Recall that the particles in BCPP either die, or make binary branching. 
To describe more general ``branching mechanism", 
we introduce a random vector 
$K=(K_x)_{x \in \zd }$ which is bounded and of 
finite range in the sense that
\bdnl{K_x}
0 \le  K_x \le b_K {\bf 1}_{\{ |x| \le r_K\}}\; \; 
\mbox{a.s. for some non-random $b_K,r_K \in [0,\8)$.}
\edn 
Let $\t^{z,i}$, ($z \in \zd$, $i \in \N^*$) 
be i.i.d. mean-one exponential random variables 
and $T^{z,i}=\t^{z,1}+...+\t^{z,i}$. 
Let also $K^{z,i}=(K_x^{z,i})_{x \in \zd }$ 
($z \in \zd$, $i \in \N^*$) be i.i.d. 
random vectors with the same distributions as $K$, 
independent of $\{\t^{z,i} \}_{z \in \zd, i \in \N^*}$. 
We suppose that the process $(\h_t)_{t \ge 0}$ 
starts from a deterministic configuration 
$\h_0=(\h_{0,x})_{x \in \zd} \in [0,\8)^{\zd}$ 
with $|\h_0|<\8$.
At time $t=T^{z,i}$, $\h_{t-}$ is replaced by $\h_t$, where 
\bdnl{h_(x,t)}
\h_{t,x}=\lef\{ \barray{ll}
K^{z,i}_0\h_{t-,z} & \mbox{if $x=z$}, \\
\h_{t-,x}+K^{z,i}_{x-z}\h_{t-,z} & 
\mbox{if $x \neq z$}.
\earray \ri.
\edn
The BCPP is a 
special case of this set-up, in which 
\bdnl{binK} 
K= \lef\{ \barray{ll}
0 & \mbox{with probability ${1 \over 2d\lm +1}$} \\
\lef(\del_{x,0}+\del_{x,e} \rig)_{x \in \zd} & 
\mbox{with probability ${\lm \over 2d\lm +1}$, for each $2d$ neighbor 
$e$ of 0.} \earray  \rig. \edn
A formal construction of the process $(\h_t)_{t \ge 0}$ can be 
given as a special case of  \cite[page 427, Theorem 1.14]{Lig85} 
via Hille-Yosida theory. In section \ref{p|ovh|=M}, 
we will also give an alternative construction of the process 
in terms of a stochastic differential equation.

We set
\bdmn
\kp_p&=&\sum_{x \in \zd}P[(K_x-\del_{x,0})^p], \; \; \; 
p=1,2,\label{kp12} \\
\ov{\h}_t &=&(\exp (-\kp_1t )\h_{t,x})_{x \in \zd}. \label{ovh_t}
\edmn
Then, 
\bdnl{|ovh|=M}
\mbox{$(|\ov{\h}_t|)_{t \ge 0}$ 
is a nonnegative martingale.}
\edn
The above martingale property can be seen by 
the same argument as in \cite[page 433, Theorem 2.2 (b)]{Lig85}. 
For the reader's convenience, we will also present a simpler 
proof in section \ref{p|ovh|=M} below.  
By (\ref{|ovh|=M}),
following limit exists almost surely:
\bdnl{ovh_8}
|\ov{\h}_\8|\st{\rm def}{=}
\lim_t|\ov{\h}_t|.
\edn
To state \Thm{CLT}, we define
\bdnl{G(x)}
G(x)=\int^\8_0P_S^0(S_t=x)dt,
\edn
where $((S_t)_{t \ge 0}, P_S^x)$ is the continuous-time random walk on 
$\zd$ starting from $x \in \zd$, with the generator 
\bdnl{L_S}
L_Sf (x)=\half \sum_{y \in \zd}\lef( P[K_{x-y}]+P[K_{y-x}] \ri)
\lef( f(y)-f(x)\rig).
\edn
As before, 
$C_{\rm b} (\rd)$ stands for the set of bounded continuous functions on $\rd$.
\Theorem{CLT}
Suppose (\ref{K_x}) and that
\bdmn
& & \mbox{the set $\{x \in \zd\; ;\; P[K_x]\neq 0\}$ contains a 
linear basis of $\rd$,} \label{K1}\\
& & \sum_{y \in \zd }P[(K_y-\del_{y,0})(K_{x+y}-\del_{x+y,0})]=0
\; \; \mbox{for all $x \in \zd \bsh \{0\}$.}
\label{K4}
\edmn
Then, referring to (\ref{kp12})--(\ref{L_S}), 
the following are equivalent:
\bds
\item[(a)] ${\kp_2 \over 2}G(0)<1$,
\item[(b)] 
${\dps \sup_{t \ge 0} P[|\ov{\h}_t|^2]<\8}$,
\item[(c)]
${\dps  \lim_{t \ra \8} \sum_{x \in \zd}
f\lef((x -mt)/\sqrt{t}\ri)\ov{\h}_{t,x}=|\ov{\h}_\8|\int_{\rd}fd\n}$  
in $\bL^2 (P)$ for all $f \in C_{\rm b} (\rd)$, 
\eds
where $ m=\sum_{x \in \zd}xP[K_x] \in \rd $ 
and $\n$ is the Gaussian measure with
\bdnl{nu}
\int_{\rd}x_id\n (x)=0, 
\; \; \; \int_{\rd}x_ix_jd\n (x)
=\sum_{x \in \zd}x_ix_jP[K_x],\; \; \; 
i,j=1,..,d.
\edn
Moreover, if ${\kp_2 \over 2}G(0)<1$, then, 
there exists $C \in (0,\8)$ such that
\bdnl{RepDec}
\sum_{x, \tl{x}\in \zd}f(x-\tl{x})P[\ov{\h}_{t,x}\ov{\h}_{t,\tl{x}}]  
\le Ct^{-d/2}|\h_0|^2\sum_{x \in \zd}f(x)
\edn
for all $t>0$ and $f :\zd \ra [0,\8)$ with $\sum_{x \in \zd}f(x)<\8$. 
\end{theorem}
The main point of \Thm{CLT} is that the condition (a), 
or equivalently (b), 
implies the central limit theorem (c) (See also \Cor{CLT} below). 
This seems to be the first result 
in which the central limit theorem for the spatial distribution 
of the particle is shown in the context of linear systems. 
Some other part of our results ((a) $\Ra$ (b), and  
\Thm{cov} below) generalizes \cite[Theorem 1]{Gri83}.
However,  this is merely a by-product 
and not a central issue in the present paper. 

The proof of \Thm{CLT}, 
which will be presented in section \ref{pCLT}, 
is roughly divided into two steps:
\bds
\item[(i)] to represent the two-point function 
$P[\h_{t,x}\h_{t,\tl{x}}]$ in terms of a continuous-time 
Markov chain on $\zd \times \zd$ via the Feynman-Kac formula 
(\Lem{FK2} and \Lem{FK3} below), 
\item[(ii)] to show the central limit theorem for the ``weighted" 
Markov chain, where the weight comes from the additive functional 
due to the Feynman-Kac formula (\Lem{perturb2} below).
\eds
The above strategy was adopted earlier by one of the authors for 
branching random walk in random 
environment \cite{Yos08}. There, the Markov chain alluded to above 
is simply the product of simple random walks on $\zd$, so that the central 
limit theorem with the Feynman-Kac weight is relatively easy.
Since the Markov chain in the present paper is no longer a random walk, 
it requires more work.
However, the good news here is that the Markov chain we have to work on is 
``close" to a random walk. In fact, we get the central limit 
theorem by perturbation from that for a random walk case. 

Some other remarks on \Thm{CLT} are in order:

\noindent {\bf 1)} 
The condition (\ref{K1}) guarantees a reasonable non-degeneracy
for the transition mechanism (\ref{h_(x,t)}). 
On the other hand, (\ref{K4}) follows from a stronger condition:
\bdnl{>K4}
P[(K_x-\del_{x,0})(K_y-\del_{y,0})]=0
\; \; \mbox{for $x,y \in \zd$ with $x \neq y$,}
\edn
which amounts to saying that 
the transition mechanism (\ref{h_(x,t)})
updates the configuration 
by ``at most one coordinate at a time". A typical examples of such $K$'s 
are given by ones which satisfy: 
$$
P(K=0)+\sum_{a \in \zd \bsh \{0\}}
P\lef( K=(\del_{x,0}+K_a\del_{x,a})_{x \in \zd}\ri)=1.
$$
These include not only BCPP but also 
models with asymmetry 
and/or long (but finite) range.

Here is an explanation for how we use the condition (\ref{K4}).
To prove \Thm{CLT}, we use a certain Markov chain 
on $\zd \times \zd$, which is introduced in \Lem{FK2} below. 
Thanks to (\ref{K4}), the Markov chain is stationary with respect 
to the counting measure on $\zd \times \zd$. 
The stationarity plays an important role in the proof of \Thm{CLT}--
see \Lem{FK3} below. \\
\noindent{\bf 2)}
Because of (\ref{K1}), the random walk 
$(S_t)$ is recurrent for 
$d=1,2$ and transient for $d \ge 3$. Therefore, 
${\kp_2 \over 2}G(0)<1$ is possible only if $d \ge 3$.
As will be explained in the proof, 
${\kp_2 \over 2}G(0)<1$ is equivalent to 
$$
P_S^0\lef[\exp \lef( {\kp_2 \over 2}\int^\8_0 \del_0 (S_t)dt\ri) \ri]<\8.
$$
\noindent{\bf 3)}
If, in particular, 
\bdnl{SRW}
P[K_x] = \lef\{ \barray{ll} 
c >0 & \mbox{for $|x|=1$,} \\  0 & \mbox{for $|x| \ge 2$,}
\earray \rig.
\edn
then, $(S_t)_{t \ge 0} \st{\mbox{\scriptsize law}}{=}
(\widehat{S}_{2dct})_{t \ge 0}$, 
where $(\widehat{S}_\cdot)$ is the simple 
random walk. Therefore, the condition (a) becomes
\bdnl{Khas2}
{\kp_2 \over 4dc(1-\pi_d)}<1.
\edn
By (\ref{binK}), the BCPP satisfies 
(\ref{K1})--(\ref{K4}). Furthermore, $\kp_2=1$ and 
we have (\ref{SRW}) with $c={\lm \over 2d\lm +1}$.  
Therefore, (\ref{Khas2}) is equivalent to (\ref{|h_8|>0}).\\
\noindent{\bf 4)}
The dual process of $(\h_t)$ above (in the sense of \cite[page 432]{Lig85})
is given by replacing the linear transform in (\ref{h_(x,t)}) by its 
transpose:
\bdnl{h_(x,t)*}
\h_{t,x}=\lef\{ \barray{ll}
\sum_{y \in \zd}K^{z,i}_{y-x}\h_{t-,y} & \mbox{if $x=z$}, \\
\h_{t-,x} & \mbox{if $x \neq z$}.
\earray \ri.
\edn 
As can be seen from the proofs, all the results in 
this paper remain true for the dual process. \\
\noindent{\bf 5)}
The central limit theorem for discrete time linear systems is 
discussed in \cite{Nak08}.

\vvs
We define the density and the replica overlap in the same way 
as (\ref{rh})--(\ref{rh^*}). Then, as 
an immediate consequence of \Thm{CLT}, we have the following 
\Corollary{CLT}
Suppose (\ref{K_x}), (\ref{K1})--(\ref{K4}) and that 
${\kp_2 \over 2}G(0)<1$. 
Then, $P [|\ov{\h}_\8|]=1$ and for all $f \in C_{\rm b} (\rd)$, 
$$
\lim_{t \ra \8} \sum_{x \in \zd}
f\lef((x -mt)/ \sqrt{t}\ri)\rh_{t,x}
=\int_{\rd}fd\n\; \; \; 
\mbox{in $P(\; \cdot \; | |\ov{\h}_\8|>0)$-probability,}
$$
where $ m=\sum_{x \in \zd}xP[K_x] \in \rd $ and 
$\n$ is the same Gaussian measure defined by (\ref{nu}).
Furthermore, 
$$
\cR_t=\cO (t^{-d/2})\; \; \; 
\mbox{as $t \nearrow \8$ in $P(\; \cdot \; | |\ov{\h}_\8|>0)$-probability}.
$$
\end{corollary}
Proof: The first statement is immediate from \Thm{CLT}(c). Taking 
$f(x)=\del_{x,0}$ in (\ref{RepDec}), we see that
$$
P[\sum_{x \in \zd}\ov{\h}_{t,x}^2 ]\le Ct^{-d/2}|\h_0|^2\; \; \; 
\mbox{for $t >0$}.
$$
This implies the second statement. 
\hfill $\Box$

\vvs
For $a \in \zd$, let $\h_t^a$ be the process starting from  
$\h_0 =(\del_{a,x})_{x \in \zd}$. 
As a by-product of \Thm{CLT}, we have the following formula 
for the covariance of $(|\ov{\h}_\8^a|)_{a \in \zd}$.
For BCPP, this formula was obtained by D. Griffeath 
\cite[Theorem 3]{Gri83}. 
\Theorem{cov}
Suppose (\ref{K_x}), (\ref{K1})--(\ref{K4}) and that 
${\kp_2 \over 2}G(0)<1$. Then,
$$
P[|\ov{\h}_\8^a||\ov{\h}_\8^b|]=1+{\kp_2G(a-b) \over 2-\kp_2G(0)}, 
\; \; \; a,b \in \zd.
$$
\end{theorem}
The proof of \Thm{cov} will be presented in section \ref{pcov}. 
We refer the reader to \cite{Yo08b} for similar formulae 
for discrete time models.
\subsection{SDE description of the process} \label{p|ovh|=M}
We now give an alternative description of the process in terms of a 
stochastic differential equation (SDE), which will be used in 
the proof of \Lem{FK2} below. 
We introduce random measures on $[0,\8) \times [0,\8)^{\zd}$ by 
\bdnl{N^z}
N^z( dsd\xi)
=\sum_{i \ge 1}{\bf 1}\{ (T^{z,i}, K^{z,i}) \in dsd\xi \}, 
\; \; \; N^z_t(dsd\xi)={\bf 1}_{\{s \le t\}}N^z(dsd\xi).
\edn
Then, $N^z$, $z \in \zd$ are independent Poisson random measures on 
$[0,\8) \times [0,\8)^{\zd}$ with the intensity 
$$
ds \times P (K \in d\xi).
$$ 
The precise definition of the process $(\h_t)_{t \ge 0}$ is then  given by 
the following stochastic differential equation:
\bdnl{sde}
\h_{t,x}=\h_{0,x}
+\sum_{z \in \zd}\int N^z_t(dsd\xi)
\lef( \xi_{x-z}-\del_{x,z}\ri)\h_{s-,z}.
\edn
By (\ref{K_x}), it is standard to see that (\ref{sde}) defines  
a unique process $\h_t=(\h_{t,x})$, ($t \ge 0$) 
and that $(\h_t)$ is Markovian. 

\vvs
\noindent {\bf Proof of (\ref{|ovh|=M})}:
Since $|\ov{\h}_t|$ is obviously nonnegative, we will prove the 
martingale property. By (\ref{sde}), we have 
$$
|\h_t|=|\h_0|
+\sum_{z \in \zd}\int N^z_t(dsd\xi)
\lef( |\xi|-1\ri)\h_{s-,z},
$$
and hence 
\bds
\item[(1)]\hspace{1cm}
${\dps |\ov{\h}_t|=|\h_0|-\kp_1\int^t_0|\ov{\h}_s|ds
+\sum_{z \in \zd}\int N^z_t(dsd\xi)
\lef( |\xi|-1\ri)\ov{\h}_{s-,z}}$.
\eds
We have on the other hand that
$$
\kp_1\int^t_0|\ov{\h}_s|ds=
\sum_{z \in \zd}\int^t_0ds \int P(K \in \xi)(|\xi|-1)\ov{\h}_{s,z}.
$$
Plugging this into (1), we see that the right-hand-side of (1) 
is a martingale.
\hfill $\Box$
\SSC{Lemmas}
\subsection{Markov chain representations for the point functions}
We assume (\ref{K_x}) throughout, but not  (\ref{K1})--(\ref{K4}) 
for the moment. 
To prove the Feynman-Kac formula for two-point function, we 
introduce some notation.

For $x,y,\tl{x},\tl{y} \in \zd$, 
\bdmn
\Gm_{x,\tl{x},y,\tl{y}}
& \st{\rm def}{=}& 
P[(K_{x-y}-\del_{x,y} )\del_{\tl{x},\tl{y}}
+(K_{\tl{x}-\tl{y}}-\del_{\tl{x},\tl{y}} )\del_{x,y}] \nn \\
& & +P[(K_{x-y}-\del_{x,y} )(K_{\tl{x}-y}-\del_{\tl{x},y}) ]\del_{y,\tl{y}},
\label{Gm2} \\
V(x) & \st{\rm def}{=}& 
\sum_{y,\tl{y} \in \zd}\Gm_{x,0,y,\tl{y}} =2\kp_1 +
\sum_{y \in \zd }P[(K_y-\del_{y,0})(K_{x+y}-\del_{x+y,0})].
\label{V(x)}
\edmn
Note that
\bdnl{V(x-x)}
V(x-\tl{x})=\sum_{y,\tl{y} \in \zd}\Gm_{x,\tl{x},y,\tl{y}}.
\edn

\vs
\noindent {\bf Remark:}
The matrix $\Gm$ introded above appears also
 in \cite[page 442, Theorem 3.1]{Lig85}, since 
it is a fundamental tool to deal with the two-point 
function of the linear system. 
However, the way we use the matrix will be different from the ones in 
the existing literature. 

\vvs
We now prove the Feynman-Kac formula for two-point 
function, which is the basis of the proof of \Thm{CLT}: 
\Lemma{FK2}
Let $(X,\tl{X})=((X_t, \tl{X}_t)_{t \ge 0}, P_{X,\tl{X}}^{x,\tl{x}})$ be 
the continuous-time Markov chain on 
$\zd \times \zd$ starting from $(x,\tl{x})$, with the generator 
$$
L_{X,\tl{X}}f (x,\tl{x})=\sum_{y,\tl{y} \in \zd}
\Gm_{x,\tl{x},y,\tl{y}} \lef( f(y,\tl{y})-f(x,\tl{x})\rig),
$$
where $\Gm_{x,\tl{x},y,\tl{y}}$ is defined by (\ref{Gm2}). 
Then, for $(t,x,\tl{x})\in [0,\8) \times \zd \times \zd$, 
\bdnl{FK2}
P[\h_{t,x}\h_{t,\tl{x}}]=P^{x,\tl{x}}_{X,\tl{X}}
\lef[ \exp \lef( \int^t_0 V(X_s-\tl{X}_s)ds\ri)
\h_{0,X_t}\h_{0,\tl{X}_t} \ri],
\edn
where $V$ is defined by (\ref{V(x)}).
\end{lemma}
Proof: 
We first show that $u(t,x,\tl{x})\st{\rm def}{=}P[\h_{t,x}\h_{t,\tl{x}}]$ 
solves the integral equation 
\bds
\item[(1)] \hspace{1cm}
${\dps u(t,x,\tl{x})-u(0,x,\tl{x})=
\int^t_0( L_{X,\tl{X}}+V(x-\tl{x}))u(s,x,\tl{x})ds.}$
\eds
By (\ref{sde}), we have 
$$
\h_{t,x}\h_{t,\tl{x}}-\h_{0,x}\h_{0,\tl{x}}
=\sum_{y \in \zd}\int N^y (dsd\xi )F_{x,\tl{x},y}(s-,\xi, \h),
$$
where
\bdnn
\lefteqn{F_{x,\tl{x},y}(s,\xi, \h)}\\
&=&
(\xi_{x-y}-\del_{x,y})\h_{s,\tl{x}}\h_{s,y}
+(\xi_{\tl{x}-y}-\del_{\tl{x},y})\h_{s,x}\h_{s,y}
+(\xi_{x-y}-\del_{x,y})(\xi_{\tl{x}-y}-\del_{\tl{x},y})\h_{s,y}^2
\ednn
Therefore, 
\bdnn
u(t,x,\tl{x})-u(0,x,\tl{x})
&=&
\sum_{y \in \zd}\int^t_0ds 
\int P[F_{x,\tl{x},y}(s,\xi, \h)]P (K \in \xi)\\
&=&
\int^t_0 \sum_{y, \tl{y}\in \zd}\Gm_{x,\tl{x},y,\tl{y}} u(s,y,\tl{y})ds\\
&\st{\scriptsize (\ref{V(x-x)})}{=}& 
\int^t_0 \lef(\sum_{y, \tl{y}\in \zd}\Gm_{x,\tl{x},y,\tl{y}} 
( u(s,y,\tl{y})-u(s,x,\tl{x}) )+V(x-\tl{x})u(s,x,\tl{x})\ri) ds \\
&=&
\int^t_0( L_{X,\tl{X}}+V(x-\tl{x}))u(s,x,\tl{x})ds.
\ednn
We next show that
\bds
\item[(2)] \hspace{1cm}
${\dps \sup_{t \in [0,T]}\sup_{x,\tl{x} \in \zd}|u(t,x,\tl{x})|<\8}$
 for any $T \in (0,\8)$.
\eds
We have by (\ref{K_x}) and (\ref{sde}) 
that, for any $p \in \N^*$, 
there exists $C_1 \in (0,\8)$ such that
$$
P[\h_{t,x}^p] \le C_1\sum_{y:|x-y| \le r_K}\int^t_0 P[\h_{s,y}^p]ds,
 \; \; \; t \ge 0.
$$
By iteration, we see that 
there exists $C_2 \in (0,\8)$ such that
$$
P[\h_{t,x}^p] \le e^{C_2t}\sum_{y \in \zd}e^{-|x-y|}(1+\h_{0,y}^p),
 \; \; \; t \ge 0,
$$
which, via Schwarz inequality,  implies (4). \\
The solution to (1) subject to (2) is unique, for each given $\h_0$. 
This can be seen by using Gronwall's inequality with respect to the 
norm $\| u\|=\sum_{x,\tl{x} \in \zd}e^{-|x|}|u(x,\tl{x})|$. 
Moreover, the RHS of (\ref{FK2}) is a solution to (1) 
subject to the bound (2). This can be seen by adapting 
the argument in \cite[page 5,Theorem 1.1]{Szn98}. Therefore, we 
get (\ref{FK2}).
\hfill $\Box$

\vvs
\noindent {\bf Remark:} The following 
Feynman-Kac formula for one-point function can be 
obtained in the same way as \Lem{FK2}:
\bdnl{FK1}
P[\h_{t,x}]=e^{\kp_1t}P^x_X[\h_{0,X_t}], \; \; (t,x)\in [0,\8) \times \zd,
\edn
where $\kp_1$ is defined by (\ref{kp12}) and 
$((X_t)_{t \ge 0}, P_X^x)$ is the continuous-time random walk on 
$\zd$ starting from $x$, with the generator
$$
L_X f (x)= \sum_{y \in \zd}P[K_{x-y}]
\lef( f(y)-f(x)\rig).
$$ 
\Lemma{stat}
We have 
\bdnl{stat}
\sum_{y,\tl{y} \in \zd}\Gm_{x,\tl{x},y,\tl{y}}
=\sum_{y,\tl{y} \in \zd}\Gm_{y,\tl{y},x,\tl{x}},
\edn
if and only if (\ref{K4}) holds. In addition, (\ref{K4}) implies that 
\bdnl{V(x)2}
V(x)=2\kp_1 +\kp_2 \del_{x,0}.
\edn
\end{lemma}
Proof: We let 
$c(x)=\sum_{y \in \zd }P[(K_y-\del_{y,0})(K_{x+y}-\del_{x+y,0})].$
Then, $c(0)=\kp_2$ and, 
\bdnn
\sum_{y,\tl{y} \in \zd}\Gm_{x,\tl{x},y,\tl{y}}
&=& 2\kp_1+c(x-\tl{x}), \; \; \; \mbox{cf. (\ref{V(x)})--(\ref{V(x-x)}),}\\
\sum_{y,\tl{y} \in \zd}\Gm_{y,\tl{y},x,\tl{x}}
&=& 2\kp_1 +\del_{x,\tl{x}}\sum_{y \in \zd }c(y).
\ednn
These imply the desired equivalence and (\ref{V(x)2}).
$\Box$

\vvs 
We assume (\ref{K4}) from here on. Then, 
by (\ref{stat}), $(X,\tl{X})$ is stationary with respect 
to the counting measure on $\zd \times \zd$. 
We denote the dual process of $(X,\tl{X})$ 
by $(Y,\tl{Y})=((Y_t,\tl{Y}_t)_{t \ge 0}, P_{Y,\tl{Y}}^{x,\tl{x}})$, 
that is, the continuous time Markov chain on 
$\zd \times \zd$ starting from $(x,\tl{x})$, with the generator 
\bdnl{L_Y*}
L_{Y,\tl{Y}}f (x,\tl{x})=\sum_{y,\tl{y} \in \zd}
\Gm_{y,\tl{y},x,\tl{x}} \lef( f(y,\tl{y})-f(x,\tl{x})\rig).
\edn
Thanks to (\ref{stat}), $L_{X,\tl{X}}$ and 
$L_{Y,\tl{Y}}$ are dual operators on $\ell^2 (\zd \times \zd)$. 

\vvs
\noindent {\bf Remark:} If we additionally 
suppose that $P[K_x^p]=P[K_{-x}^p]$ for $p=1,2$ and $x \in \zd$, 
then, $\Gm_{x,\tl{x},y,\tl{y}}=\Gm_{y,\tl{y},x,\tl{x}}$ for 
all $x,\tl{x},y,\tl{y} \in \zd$. Thus, $(X,\tl{X})$ and 
$(Y,\tl{Y})$ are the same in this case.

\vvs
The relative motion $Y_t-\tl{Y}_t$ of the components of 
$(Y,\tl{Y})$ is nicely identified by: 
\Lemma{Y-Y}
$((Y_t-\tl{Y}_t)_{t \ge 0}, P_{Y,\tl{Y}}^{x,\tl{x}})$ and 
$((S_{2t})_{t \ge 0}, P_S^{x-\tl{x}})$ (cf. (\ref{L_S})) have the same law.
\end{lemma}
Proof: 
Since $(Y,\tl{Y})$ is shift invariant, in the sense that 
$\Gm_{x+v,\tl{x}+v,y+v,\tl{y}+v} =\Gm_{x,\tl{x},y,\tl{y}} $ 
for all $v \in \zd$, $((Y_t-\tl{Y}_t)_{t \ge 0}, P_{Y,\tl{Y}}^{x,\tl{x}})$ 
is a Markov chain. Moreover, its 
jump rate is computed as follows.
For $x \neq y$, 
\bdmn
\sum_{z \in \zd}\Gm_{y+z,z,x,0}
& = & P[K_{x-y}]+P[K_{y-x}] +\del_{x,0}\sum_{z \in \zd}
P[(K_{y+z}-\del_{y,z})(K_{z}-\del_{0,z})] \nn \\
& \st{\mbox{\scriptsize (\ref{K4})}}{=} & P[K_{x-y}]+P[K_{y-x}]. \nn
\edmn
\hfill $\Box$

\vvs
To prove \Thm{CLT}, the use of \Lem{FK2} is made not in itself, but 
via the following lemma. It is the proof of this lemma, where the 
duality of $(X,\tl{X})$ and $(Y,\tl{Y})$ plays its role.
\Lemma{FK3}
For a bounded $g:\zd \times \zd \ra \R$, 
\bdmn
\lefteqn{\sum_{x,\tl{x} \in \zd}
P[\ov{\h}_{t,x}\ov{\h}_{t,\tl{x}}]g(x,\tl{x})} \nn \\
&=&\sum_{x,\tl{x} \in \zd}\h_{0,x}\h_{0,\tl{x}}
P_{Y,\tl{Y}}^{x,\tl{x}}
\lef[ \exp \lef( \kp_2\int^t_0 \del_0(Y_s-\tl{Y}_s)ds\ri)
g(Y_t,\tl{Y}_t) \ri].\label{ell_2}
\edmn
In particular, for a bounded $f:\zd  \ra \R$, 
\bdnl{FK3}
\sum_{x,\tl{x} \in \zd}P[\ov{\h}_{t,x}\ov{\h}_{t,\tl{x}}]f(x-\tl{x})
=\sum_{x,\tl{x} \in \zd}\h_{0,x}\h_{0,\tl{x}}
P_S^{x-\tl{x}}\lef[ \exp \lef( {\kp_2 \over 2}\int^{2t}_0 \del_0(S_u)du\ri) 
f(S_{2t})\ri].
\edn
\end{lemma}
Proof: 
It follows from \Lem{FK2} and (\ref{V(x)2}) that
\bds
\item[(1)]
\hspace{1cm} ${\dps 
\mbox{LHS of (\ref{ell_2})}=\sum_{x,\tl{x} \in \zd}
P_{X,\tl{X}}^{x,\tl{x}}
\lef[ \exp \lef( \kp_2\int^t_0 \del_0(X_s-\tl{X}_s)ds\ri)
\h_{0,X_t}\h_{0,\tl{X}_t} \ri]g(x,\tl{x}).}$
\eds
We now observe that the operators 
\bdnn
f (x,\tl{x}) & \mapsto &
P_{X,\tl{X}}^{x,\tl{x}}
\lef[ \exp \lef( \kp_2\int^t_0 \del_0(X_s-\tl{X}_s)ds\ri)
f(X_t, \tl{X}_t) \ri], \\
f (x,\tl{x}) & \mapsto &
P_{Y,\tl{Y}}^{x,\tl{x}}
\lef[ \exp \lef( \kp_2\int^t_0 \del_0(Y_s-\tl{Y}_s)ds\ri)
f(Y_t, \tl{Y}_t) \ri]
\ednn
are dual to each other with respect to the counting measure on
 $\zd \times \zd $. Therefore, 
$$
\mbox{RHS of (1)}=\mbox{RHS of (\ref{ell_2})}.
$$
Taking $g (x,\tl{x})=f(x-\tl{x})$ in particular, we have by (\ref{ell_2}) 
and \Lem{Y-Y} that
\bdnn
\mbox{LHS of (\ref{FK3})}
&=&\sum_{x,\tl{x} \in \zd}\h_{0,x}\h_{0,\tl{x}}P^{x,\tl{x}}_{Y,\tl{Y}}
\lef[ \exp \lef( \kp_2\int^t_0 \del_0(Y_s-\tl{Y}_s)ds\ri)
f(Y_t-\tl{Y}_t)\ri] \\
&=& \sum_{x,\tl{x} \in \zd}\h_{0,x}\h_{0,\tl{x}}
P_S^{x-\tl{x}}\lef[ \exp \lef( \kp_2\int^t_0 \del_0(S_{2u})du\ri)
f(S_{2t}) \ri]
=\mbox{RHS of (\ref{FK3})}.
\ednn
\hfill $\Box$

\vs
\noindent{\bf Remark:}
In the case of BCPP, D. Griffeath obtained 
a Feynman-Kac formula for 
$$
\sum_{y \in \zd}P[\h_{t,x}\h_{t,\tl{x}+y}]
$$
\cite[proof of Theorem 1]{Gri83}. However, this does not 
seem to be enough for our purpose. 
Note that the Feynman-Kac formulae 
in the present paper (\Lem{FK2} and \Lem{FK3}) are stronger, 
since they give the expression for each summand of the above summation. 
\subsection{Central limit theorems for Markov chains}
We prepare central limit theorems for Markov chains, which is obtained 
by perturbation of random walks. 
\Lemma{perturb1}
Let $((Z_t)_{t \ge 0}, P^x)$ be a continuous-time random walk 
on $\zd$ starting from $x$, 
with the generator 
$$
L_Z f(x)=\sum_{y \in \zd}a_{y-x}(f(y)-f(x)),
$$ 
where we assume that
$$
\sum_{x \in \zd}|x|^2a_x<\8.
$$
Then, for any $B \in \s [Z_u\; ; \; u \in [0,\8)]$, 
$x \in \zd$, and $f \in C_{\rm b}(\rd)$,
$$
\lim_{t \ra \8}P^x [f((Z_t -mt)/\sqrt{t}):B]
=P^x (B)\int_{\rd}fd\n,
$$
where $m=\sum_{x \in \zd}xa_x$ and $\n$ is the Gaussian measure with
\bdnl{nua}
\int_{\rd}x_id\n (x)=0, 
\; \; \; \int_{\rd}x_ix_jd\n (x)
=\sum_{x \in \zd}x_ix_ja_x,\; \; \; 
i,j=1,..,d.
\edn
\end{lemma}
Proof: 
By subtracting a constant, we may assume that $\int_{\rd}fd\n =0$. 
We first consider the case that 
$B \in \cF_s \st{\rm def}{=}\s [Z_u\; ; \; u \in [0,s]]$ 
for some $s \in (0,\8)$. It is easy to see 
from the central limit theorem for $(Z_t)$ that for any  $x \in \zd$,
$$
\lim_{t \ra \8}P^x [f((Z_{t-s} -mt)/\sqrt{t})]=0.
$$
With this and the bounded convergence theorem, we have 
$$
P^x[f((Z_t -mt)/\sqrt{t}):B]
=P^x[P^{Z_s}[f((Z_{t-s} -mt)/\sqrt{t})]:B]
\lra 0\; \; \mbox{as $t \nearrow \8$}.
$$
Next, we take $B \in \s [Z_u\; ; \; u \in [0,\8)]$. 
For any $\e>0$, there exist $s \in (0,\8)$ and $\tl{B} \in \cF_s$  
such that $P^x[|{\bf 1}_B-{\bf 1}_{\tl{B}}|] <\e$. 
Then, by what we already have seen, 
$$
\suplim_{t \ra \8}
P^x[f((Z_t -mt)/\sqrt{t}):B] 
 \le \suplim_{t \ra \8}P^x[f((Z_t -mt)/\sqrt{t}):\tl{B}] 
+\| f \|\e  = \| f \|\e,
$$
where $\| f\|$ is the sup norm of $f$. Similarly, 
$$
\inflim_{t \ra \8}
P^x[f((Z_t -mt)/\sqrt{t}):B] \ge   -\| f \|\e.
$$
Since $\e>0$ is arbitrary, we are done.
\hfill $\Box$
\Lemma{perturb2}
Let $Z=((Z_t)_{t \ge 0}, P^x)$ be as in \Lem{perturb1} and 
and $D \sub \zd$ be transient for $Z$.
On the other hand, let $\tl{Z}=((\tl{Z}_t)_{t \ge 0}, \tl{P}^x)$ 
be the continuous-time Markov chain on $\zd$ starting from $x$, 
with the generator 
$$
L_{\tl{Z}} f(x)=\sum_{y \in \zd}\tl{a}_{x,y}(f(y)-f(x)),
$$ 
where we assume that $\tl{a}_{x,y}=a_{y-x}$ if $x \not\in D \cup \{y\}$ 
and that $D$ is also transient for $\tl{Z}$. 
Furthermore, we assume that 
a function $v:\zd \ra \R$ satisfies
\bdnn
& & \mbox{$v \equiv 0$ outside $D$,}\\
& & \tl{P}^z\lef[\exp \lef( \int_0^\8 |v(\tl{Z}_t)|dt\ri)\ri]<\8\; \;
\mbox{for some $z \in \zd$}.
\ednn
Then, for $f \in C_{\rm b}(\rd)$, 
$$
\lim_{t \ra \8}\tl{P}^z
\lef[\exp \lef( \int_0^t v(\tl{Z}_u)du\ri)f((\tl{Z}_t-mt)/\sqrt{t})\ri]
=\tl{P}^z\lef[\exp \lef( \int_0^\8 v(\tl{Z}_t)dt\ri)\ri]\int_{\rd}fd\n,
$$
where $\n$ is the Gaussian measure such that (\ref{nua}) holds.
\end{lemma}
Proof:
Define  
\bdnn
H_D(\tl{Z})&=&\inf \{ t \ge 0\; ; \; \tl{Z}_t \in D\}, \; \; 
T_D(\tl{Z})=\sup \{ t \ge 0\; ; \; \tl{Z}_t \in D\},\\
e_t&=&\exp \lef( \int_0^t v(\tl{Z}_s)ds \ri).
\ednn
Then, for $s<t$, 
\bdmn
\lefteqn{\tl{P}^z\lef[e_tf((\tl{Z}_t-mt)/\sqrt{t})\ri]} \nn \\
&=&\tl{P}^z\lef[e_tf((\tl{Z}_t-mt)/\sqrt{t}):T_D(\tl{Z})<s\ri]+\e_{s,t}\nn\\
&=&\tl{P}^z\lef[e_sf((\tl{Z}_t-mt)/\sqrt{t}):T_D(\tl{Z})<s\ri] +\e_{s,t} \nn\\
&=& \tl{P}^z\lef[e_s{\bf 1}_{\tl{Z}_s \not\in D}
\tl{P}^{\tl{Z}_s}\lef[f((\tl{Z}_{t-s}-mt)/\sqrt{t}):H_D(\tl{Z})=\8 \ri]\ri]
+\e_{s,t},\label{est}
\edmn
where
\bdnn
|\e_{s,t}| 
& = & 
\lef| \tl{P}^z\lef[e_tf((\tl{Z}_t-mt)/\sqrt{t}):T_D(\tl{Z})\ge s\ri]\rig| \\
& \le &\| f \|\tl{P}^z\lef[
\exp \lef( \int_0^\8 |v(\tl{Z}_t)|dt\ri):T_D(\tl{Z})\ge s\ri] \ra 0\; \; 
\mbox{as $s \ra \8$.}
\ednn
We now observe that
$$
\tl{P}^x( \; \cdot \; | H_D(\tl{Z})=\8)=P^x( \; \cdot \; | H_D(Z)=\8)\; \; \; 
\mbox{for $x \not\in D$,}
$$
where $H_D(Z)$ is defined similarly as $H_D(\tl{Z})$.
Hence, for $x \not\in D$ and fixed $s>0$, we have by \Lem{perturb1} that
$$
\lim_{t \ra \8}\tl{P}^x\lef[f((\tl{Z}_{t-s}-mt)/\sqrt{t}):H_D(\tl{Z})=\8 \ri]
= \tl{P}^x[H_D(\tl{Z})=\8] \int_{\rd}fd\n.
$$
Therefore, 
\bdnn
\lefteqn{\lim_{t \ra \8}
\tl{P}^z\lef[e_s{\bf 1}_{\tl{Z}_s \not\in D}
\tl{P}^{\tl{Z}_s}\lef[f((\tl{Z}_{t-s}-mt)/\sqrt{t}):H_D(\tl{Z})=\8 \ri]\ri]}\\
&=&\tl{P}^z\lef[e_s{\bf 1}_{\tl{Z}_s \not\in D}
\tl{P}^{\tl{Z}_s}[H_D(\tl{Z})=\8] \rig]\int_{\rd}fd\n \\
&=& \tl{P}^z\lef[e_s :T_D (\tl{Z})<s \rig]\int_{\rd}fd\n. \\
\ednn
Thus, letting $t \ra \8$ first, and then $s \ra \8$, in (\ref{est}), 
we get the lemma.
\hfill $\Box$
\subsection{A Nash type upper bound for the Schr\"odinger semi-group}
We will use the following lemma to prove (\ref{RepDec}). The lemma 
can be generalized to symmetric Markov chains on more general graphs.
However, we restrict ourselves to random walks on $\zd$, since 
it is enough for our purpose.
\Lemma{Nash}
Let $((Z_t)_{t \ge 0}, P^x)$ be continuous-time random walk 
on $\zd$
with the generator:
$$
L_Z f(x)=\sum_{y \in \zd}a_{y-x}(f(y)-f(x)),
$$ 
where we assume that
\bdnn
& &
\mbox{the set $\{x \in \zd\; ;\; a_x\neq 0\}$ is bounded and contains a 
linear basis of $\rd$,}\\
& & a_x=a_{-x}\; \; \mbox{for all $x \in \zd$},
\ednn
Let $v:\zd \ra \R$ be a function such that
$$
C_v \st{\rm def}{=}
\sup_{x \in \zd}P^x\lef[\exp \lef( \int_0^\8 |v(Z_t)|dt\ri)\ri]<\8.
$$
Then, there exists $C \in (0,\8)$ such that
\bdnl{RepDec2}
\sup_{x \in \zd} P^x\lef[\exp \lef( \int_0^t v(Z_u)du\ri)f(Z_t)\ri]
\le Ct^{-d/2}\sum_{x \in \zd}f(x)
\edn
for all $t>0$ and $f :\zd \ra [0,\8)$ with $\sum_{x \in \zd}f(x)<\8$. 
\end{lemma}
Proof: We adapt the argument in \cite[Lemma 3.1.3]{CY04}.
For a bounded function $f:\zd \ra \R$, we introduce
\bdnn
(T_tf)(x)&=&P^x \lef[ 
\exp \lef( \int^t_0v(Z_u)du\rig)f(Z_t)\rig],
\; \; \; x \in \zd,\\
T_t^hf&=&\frac{1}{h}T_t [fh], \; \; 
\mbox{where ${\dps h(x)=P^x\lef[\exp \lef( \int_0^\8 v(Z_t)dt\ri)\ri]}$.}
\ednn
Then,
$(T_t)_{t \ge 0}$ extends to a symmetric,
strongly continuous semi-group on $\ell^2(\zd)$.
We now consider the measure $\sum_{x \in \zd}h(x)^2\del_x$ on $\zd$, 
and denote by $(\ell^{p,h}(\zd), \| \; \cdot\;  \|_{p,h})$ the associated 
$\bL^p$-space. Then, it is standard 
(e.g., proofs of \cite[page 74, Theorem 3.10]{ChZh95} and 
\cite[page 16, Proposition 3.3]{Szn98}) 
to see that $(T_t^h)_{t \ge 0}$
defines a symmetric strongly continuous semi-group on $\ell^{2,h}(\zd)$ 
and that for $f \in \ell^{2,h}(\zd)$,
\bdnn
\cE^h (f,f)  & \st{\rm def.}{=} & \lim_{t \searrow 0}{1 \over t}
\sum_{x\in \zd} f (x)(f-T_t^hf)(x)h(x)^2  \\
&=& \half \sum_{x,y \in \zd}a_{y-x}| f(y)-f(x)|^2h(x)h(y).
\ednn
By the assumptions on $(a_x)$, we have the Sobolev inequality:
\bds
\item[(1)] \hspace{1cm}
${\dps \sum_{x \in \zd}| f (x)|^{\frac{2d}{d-2}} \le c_1
\lef( \half\sum_{x,y \in \zd}a_{y-x}| f(y)-f(x)|^2 \rig)^{\frac{d}{d-2}}
}$ for all $f \in \ell^2 (\zd)$, 
\eds
where $c_1\in (0,\8)$ is independent of $f$. 
This can be seen via an isoperimetric inequality \cite[page 40, (4.3)]{Woe00}. 
We have on the other hand that
\bds
\item[(2)] \hspace{1cm}
$ 1/C_v\le h(x) \le C_v.$
\eds  
We see from (1) and (2) that 
$$
\sum_{x \in \zd}| f (x)|^{\frac{2d}{d-2}}h(x)^2\le c_2
\cE^h (f,f) ^{\frac{d}{d-2}}
\; \; \; \mbox{for all $f \in \ell^{2,h} (\zd)$,}
$$
where $c_2\in (0,\8)$ is independent of $f$.
This implies that
there is a constant $C$ such that
$$
\| T^h_t \|_{2 \ra \8,h} \le Ct^{-d/4}
\; \; \; \mbox{for all $t >0$},
$$
e.g.,\cite[page 75, Theorem 2.4.2]{Dav89},
where $\| \cdot \|_{p \ra q,h}$ denotes
the operator norm from $\ell^{p,h}(\zd)$ to $\ell^{q,h}(\zd)$.
Note that $\| T^h_t   \|_{1 \ra 2,h}
=\| T^h_t   \|_{2 \ra \8,h}$
by duality. We therefore have via semi-group property that
\bds
\item[(3)] \hspace{1cm}
$\| T^h_t \|_{1 \ra \8,h} \le
\| T^h_{t/2} \|_{2 \ra \8,h}^2 \le C^2t^{-d/2}$ for all $t >0$.
\eds
Since $T_tf=h T^h_t[f/h]$,
the desired bound (\ref{RepDec2})
follows from (2) and (3).
\hfill $\Box$
\SSC{Proof of \Thm{CLT} and \Thm{cov}}
\subsection{Proof of \Thm{CLT}} \label{pCLT}
(a) $\LRa$ (b): 
Define 
$$
h (x)=P_S^x\lef[\exp \lef( {\kp_2 \over 2}\int^\8_0 \del_0 (S_t)dt\ri)\ri].
$$
Since $\max_{x \in \zd}h(x)=h(0)$, we have that
$$
\sup_{t \ge 0} P[|\ov{\h}_t|^2]\st{\scriptstyle (\ref{FK3})}{=}
\sum_{x,\tl{x}\in \zd}\h_{0,x}\h_{0,\tl{x}}h(x-\tl{x})
\lef\{ \barray{l} 
\le h(0)|\h_0|^2, \\
\ge h(0)\sum_{x \in \zd}\h_{0,x}^2.
\earray \rig.
$$
Therefore, it is enough to show that (a) is equivalent to $h(0)<\8$.
In fact, Khas'minskii's lemma 
(e.g., \cite[page 71]{ChZh95} or \cite[page 8]{Szn98}) says that 
$\sup_{x \in \zd}h(x)<\8$ if  ${\kp_2 \over 2}\sup_{x \in \zd}G(x)<1$. 
Since $\max_{x \in \zd}G(x)=G(0)$, (a) implies that $h(0)<\8$.
On the other hand, we have that
$$
\exp \lef( {\kp_2 \over 2}\int^t_0 \del_0 (S_s)ds\rig)
=1+{\kp_2 \over 2}\int^t_0\del_{0,S_s}
\exp \lef( {\kp_2 \over 2}\int^t_s \del_0 (S_u)du\rig)ds, 
$$
and hence that $h (x)=1+{\kp_2 \over 2} h(0)G(x)$. Thus, $h(0)<\8$ 
implies (a) and that 
\bdnl{h(x)}
h(x)=1+{\kp_2G(x) \over 2-\kp_2G(0)}.
\edn
(a),(b) $\Ra$ (c): 
Since (b) implies that 
$\lim_{t \ra \8}|\ov{\h}_t|=|\ov{\h}_\8|$ in $\bL^2 (P)$,
 it is enough to prove that
$$
U_t\st{\rm def.}{=}\sum_{x \in \zd}\ov{\h}_{t,x}f \lef( (x-mt)/\sqrt{t}\rig)
\lra 0\; \; \; \mbox{in $\bL^2 (P)$ as $t \nearrow \8$}
$$
for $f \in C_{\rm b} (\rd)$ such that  $\int_{\rd}fd\n =0$.

We set $f_t(x, \tl{x}) =f ((x-m)/\sqrt{t})f ((\tl{x}-m)/\sqrt{t})$. 
Then, by \Lem{FK3},
$$
P[U_t^2]
 =  \sum_{x,\tl{x} \in \zd}
P[\ov{\h}_{t,x}\ov{\h}_{t,\tl{x}}]f_t(x, \tl{x})
= \sum_{x,\tl{x} \in \zd}\h_{0,x}\h_{0,\tl{x}}
P_{Y,\tl{Y}}^{x,\tl{x}}\lef[ e_tf_t(Y_t, \tl{Y}_t) \ri],
$$
where $e_t=\exp \lef( \kp_2\int^t_0 \del_0(Y_s-\tl{Y}_s)ds\ri)$.
Note that by \Lem{Y-Y} and (a), 
\bds
\item[(1)]\hspace{1cm}
${\dps P_{Y,\tl{Y}}^{x,\tl{x}}\lef[ e_\8\ri] =h(x-\tl{x}) \le h(0) <\8}$.
\eds
Since $|\h_0|<\8$, it is enough to prove that for each $x,\tl{x} \in \zd$
$$
\lim_{t \ra \8}
P_{Y,\tl{Y}}^{x,\tl{x}}\lef[ e_tf_t(Y_t, \tl{Y}_t) \ri]=0.
$$
To prove this, we apply \Lem{perturb2} to the Markov chain 
$\tl{Z}_t\st{\rm def.}{=}(Y_t,\tl{Y}_t)$ and the random walk 
$(Z_t)$ on $\zd \times \zd$ with the generator 
$$
L_Zf (x,\tl{x})
=\sum_{y,\tl{y} \in \zd}a_{x,\tl{x},y,\tl{y}} 
\lef( f(y,\tl{y})-f(x,\tl{x})\rig)\; \; \; 
\mbox{with}\; \; 
a_{x,\tl{x},y,\tl{y}}=\lef\{\barray{ll}
P[K_{\tl{y}-\tl{x}}]
& \mbox{if $x=y$ and $\tl{x} \neq \tl{y}$,}\\
P[K_{y-x}] 
& \mbox{if $x \neq y$ and $\tl{x} =\tl{y}$,}\\
0 & \mbox{if otherwise}.\earray \rig.
$$
Let $D=\{ (x,\tl{x}) \in \zd \times \zd\; ; \; x=\tl{x}\}$. 
Then, 
\bds
\item[(2)] $a_{x,\tl{x},y,\tl{y}}=\Gm_{y,\tl{y},x,\tl{x}}$ if 
$ (x,\tl{x}) \not\in D \cup \{ (y,\tl{y}) \}$, 
\eds
since
\bdnn
\lefteqn{\Gm_{y,\tl{y},x,\tl{x}} } \\
& =& 
P[(K_{y-x}-\del_{y,x} )\del_{\tl{y},\tl{x}}
+(K_{\tl{y}-\tl{x}}-\del_{\tl{y},\tl{x}} )\del_{y,x}
+(K_{y-x}-\del_{y,x} )(K_{\tl{y}-x}-\del_{\tl{y},x})\del_{x,\tl{x}} ]. \nn 
\ednn
Moreover, by (\ref{K1}), 
\bds
\item[(3)] $D$ is transient both for $(Z_t)$ and for $(\tl{Z}_t)$.
\eds
 Finally, the Gaussian measure $\n \otimes \n $ 
is the limit law in the central limit theorem 
for the random walk $(Z_t)$.
Therefore, by (1)--(3) and \Lem{perturb2},
$$
\lim_{t \ra \8}
P_{Y,\tl{Y}}^{x,\tl{x}}\lef[ e_tf_t(Y_t, \tl{Y}_t)  \ri]
=P_{Y,\tl{Y}}^{x,\tl{x}}\lef[ e_\8\ri]
\lef( \int_{\rd}fd\n \rig)^2=0.
$$
(c) $\Ra$ (b): This can be seen by taking $f \equiv 1$.\\
(\ref{RepDec}):By (\ref{FK3}),
$$
\sum_{x,\tl{x} \in \zd}P[\ov{\h}_{t,x}\ov{\h}_{t,\tl{x}}]f(x-\tl{x})
=\sum_{x,\tl{x} \in \zd}\h_{0,x}\h_{0,\tl{x}}P^{x-\tl{x}}_S\lef[ 
\exp \lef( {\kp_2 \over 2}\int^{2t}_0 \del_0(S_u)du\ri) f(S_{2t})\ri].
$$
We apply \Lem{Nash} to the right-hand-side to get  (\ref{RepDec}).
\hfill $\Box$
\subsection{Proof of \Thm{cov}} \label{pcov}
By the shift-invariance, we may assume that $b=0$. We have by \Lem{FK3} 
that
$$
P[\ov{\h}^a_{t,x}\ov{\h}^0_{t,\tl{x}}]
=P_{Y,\tl{Y}}^{a,0}
\lef[ \exp \lef( \kp_2 \int^{t}_0 \del_0(Y_u-\tl{Y}_u)du\ri) 
:(Y_t,\tl{Y}_t)=(x,\tl{x}) \ri],
$$
and hence by \Lem{Y-Y} that
$$
P[|\ov{\h}_t^a||\ov{\h}_t^0|]
=P^{a,0}_{Y,\tl{Y}}
\lef[ \exp \lef( \kp_2 \int^{2t}_0 \del_0(Y_u-\tl{Y}_u)du\ri)\ri]
=P^{a}_S\lef[ \exp \lef( {\kp_2 \over 2}\int^{2t}_0 \del_0(S_u)du\ri)\ri].
$$
By \Thm{CLT}, both $|\ov{\h}_t^a|$ and 
$|\ov{\h}_t^0|$ are convergent in $\bL^2 (P)$ if $ {\kp_2 \over 2}G(0)<1$.
Therefore, letting $t \nearrow \8$, we conclude that
$$
P[|\ov{\h}_\8^a||\ov{\h}_\8^0|]=
P^{a}_S\lef[ \exp \lef( {\kp_2 \over 2}\int^{\8}_0 \del_0(S_u)du\ri)\ri]
\st{\mbox{\scriptsize (\ref{h(x)}) }}{=}1+{\kp_2G(a) \over 2-\kp_2G(0)}.
$$
\hfill $\Box$

\small

\vvs
\noindent{\bf Acknowledgements:}
The authors thank Shinzo Watanabe for the improvement of \Lem{perturb1}.


\end{document}